\documentclass[12pt]{amsart}
\usepackage{amssymb,amsfonts,lineno,amsthm,amscd,stmaryrd}
\usepackage[leqno]{amsmath}
\usepackage[numbers]{natbib}

\usepackage[body={16cm,23.5cm}, top=3cm,left=2.5cm,right=2.5cm,bottom=3cm]{geometry}
\theoremstyle{plain}
\newtheorem{thm}{Theorem}

\newtheorem{lem}[thm]{Lemma}
\newtheorem{prop}[thm]{Proposition}
\theoremstyle{definition}

\newtheorem{remark}[thm]{Remark}

\newcommand{\Puccisub}{{\mathcal{P}}^+_{\lambda,\Lambda}}

\newcommand{\ep}{\varepsilon}
\newcommand{\R}{\ensuremath{\mathbb{R}}}
\newcommand{\rn}{\ensuremath{\mathbb{R}^N}}

\begin{document}
\title[A new approach to Liouville theorems]{A new approach to Liouville theorems for elliptic inequalities}
\author{Scott N. Armstrong}
\address{Department of Mathematics\\ The University of Chicago\\ 5734 S. University Avenue
Chicago, Illinois 60637.}
\email{armstrong@math.uchicago.edu}
\author{Boyan Sirakov}
\address{UFR SEGMI, Universit\'e Paris 10\\
92001 Nanterre Cedex, France \\
and CAMS, EHESS \\
54 bd Raspail \\
75270 Paris Cedex 06, France}
\email{sirakov@ehess.fr}
\date{\today}

\maketitle

\section{Introduction}

In this article we review a new method for proving the nonexistence of positive solutions of elliptic inequalities in unbounded domains in $\rn$, $N\geq 2$, which was recently introduced by the authors \cite{AS}. For clarity and ease of exposition here we consider only domains which contain the infinite point of $\rn$, that is, domains in the form $\R^N\! \setminus\! G$, where $G$ is an arbitrary bounded set. To further simplify the presentation we are going to state our results for the following simple but important and widely studied inequalities
\begin{equation} \label{peq}
-\Delta_p u \geq f(u) \quad \mbox{in} \ \R^N\! \setminus\! G,
\end{equation}
and
\begin{equation} \label{neq}
-\Puccisub (D^2 u) \geq f(u) \quad \mbox{in} \ \R^N\! \setminus\! G,
\end{equation}
where  $f : (0,\infty) \to (0,\infty)$ is some continuous map.  Our results, which are new even for the standard semilinear inequality $-\Delta u \ge f(u)$, provide sharp hypotheses on $f$ under which \eqref{peq} and \eqref{neq} have no positive solutions.

We have taken \eqref{peq} as an example of an inequality in divergence form, whose weak solutions are naturally defined in Sobolev spaces. Here $-\Delta_p$ denotes the $p$-Laplacian and $p> 1$.

In \eqref{neq} $\Puccisub$ denotes the Pucci maximal operator $\Puccisub (M)= \sup_{\lambda I\le A\le \Lambda I}\mathrm{tr}(-AM)$, for some positive constants $\lambda\le \Lambda$. This equation is in non-divergence form, and its weak solutions are naturally defined in the viscosity sense. Note the non-existence of solutions of \eqref{neq} implies any semi-linear inequality $-a_{ij}(x)\partial_{ij}u\ge f(u)$ has no solutions either, provided the eigenvalues of the matrix $(a_{ij}(x))$ lie in the interval $[\lambda, \Lambda]$.

  We first state the result we obtain on \eqref{peq}. The type of condition we impose on $f$ varies depending on how $p$ compares to $N$.

\begin{thm}[\cite{AS}] \label{lio}
Denote
\begin{equation*}
p_* : = \frac{N(p-1)}{N-p} \quad \mbox{if} \ p \neq N.
\end{equation*}
Assume that $f:(0,\infty) \to (0,\infty)$ is continuous and satisfies
\begin{enumerate}
\item[(i)]  if $p < N$, then $\liminf_{t \to 0}\, t^{-p_*}f(t) > 0$;
\item[(ii)] if $p=N$, then $\liminf_{t \to \infty}\, e^{a t} f(t) > 0$ for each $a> 0$; and
\item[(iii)] if $p> N$, then $\liminf_{t \to \infty} \, t^{|p_*|}f(t) > 0$.
\end{enumerate}
Then the inequality \eqref{peq} has no positive weak solution.
\end{thm}

This theorem is sharp: for instance the model inequalities $-\Delta_p u \geq u^{p_*+\ep}$ have positive solutions in every exterior domain if $p<N$ and $\ep>0$, or if $p>N$ and $\ep<0$ (resp.  $-\Delta_p u \geq e^{-au}$ has solutions if $p=N$, for each $a>0$; see the end of the paper).

The study of the nonexistence of positive supersolutions of elliptic equations and systems has a rich literature. While we do not give extensive references here, referring instead to the more complete bibliography in our paper \cite{AS} as well as to \cite{V,MP,SZ,KLS,AM}, we do mention that special cases of Theorem~\ref{lio} have been proved among other things by Gidas \cite{G},  Ni and Serrin \cite{NS}, Bidaut-Veron \cite{B}, Bidaut-Veron and Pohozaev \cite{BP}, Serrin and Zou \cite{SZ} and more recently by d'Ambrosio and Mitidieri  \cite{AM}. The previous methods for proving Liouville-type results like Theorem~\ref{lio} have involved either assembling  delicate integral identities using the integral formulation of the equation or, should the symmetries of the equation permit, "radializing" the equation, that is, showing that the spherical mean of an eventual solution satisfies an ODE without solutions. Our technique, which will be developed below, is rather different and relies on some simple ideas related to the maximum principle.

\smallskip

What is striking about Theorem~\ref{lio} at first glance is how little is required of the function~$f$. Only local conditions are imposed on the behavior of $f$, in the sense that we demand only that $f(t)$ either grow fast enough near $t=0$ or decay slowly enough near $t=\infty$, but allow arbitrary behavior elsewhere. In constrast, most of the previous papers considered the case $f(t) = t^q$, $q>0$. To our knowledge, only the hypothesis (i) in Theorem \ref{lio} has appeared before, for the first time in \cite{NS} for decaying solutions, and recently in \cite{AM} for differential inequalities holding in the whole space $\R^N$. The possibility of allowing nonlinearities which decay at infinity in the case $p\ge N$ has not been observed (except for the trivial case of an inequality in $\rn$ where $p$-superharmonic functions do not exist).

\smallskip

The conditions (i) - (iii) can be best explained if we remember the dilative scaling of the equation in the model case $f(t) =t^q$ for $q\not=p-1$. As it is easy to check, if $u$ is a solution of
\begin{equation} \label{mod}
-\Delta_p u \ge u^q \quad \mbox{in} \ \R^N \setminus B_1,
\end{equation}
then for any $s > 0$ the rescaled function $u_s(x):= s^{q^*} u(sx)$ is a supersolution of the same equation in the domain $\R^N \! \setminus \!B_{1/s}$, provided we set the scaling exponent to be $$q^* : = p/(q-p+1).$$ The question of existence or nonexistence of positive supersolutions of \eqref{mod} turns out to depend on the  competition between this scaling exponent $q^*$ and the homogeneity $$\alpha^*= (N-p)/(p-1)$$ of the fundamental solution $\Phi=\Phi_p(x)$ of the $p$-Laplace equation, which we recall is given by
\begin{equation*}
\Phi_p(x)  =
\pm|x|^{-\alpha^*} \; \mbox{ if }\; \ \alpha^*\not=0, \quad
\Phi_p(x)  = \pm\log |x| \; \mbox{ if }\;  \ \alpha^*=0.
\end{equation*}
For example, for \eqref{mod} if $q>p-1$ and $p <N$ condition (i) is equivalent to the inequality $0<\alpha^*\le q^*$. Similarly if $q<p-1$ and $p >N$ condition (iii) requires that $q^*\le \alpha^*<0$. This point of view also explains why the conditions in Theorem~\ref{lio} are sharp: to find a supersolution (e.g., in the model case $f(t) = t^q$) one needs only to slightly modify the fundamental solution $\Phi_p$ by bending it in an appropriate way. A first discussion on  the interplay between $\alpha^*$ and $q^*$ appeared in our earlier paper \cite{AS0}, where we used an argument based on a "linearization" to study some particular inequalities in non-divergence form.

\smallskip

Let us now state the result we obtain on the inequality \eqref{neq}. Dividing the inequality by $\Lambda$ we can assume $\Lambda=1$. We will also assume we are in the non-trivial case $\lambda<1$ (the case $\lambda=1$ is covered by Theorem \ref{lio} with $p=2$).

Observe that a nonexistence result for \eqref{neq} implies the rather strong assertion that all semilinear inequalities with fixed ellipticity constants and $L^\infty$-bounds for the coefficients have no solutions at infinity. So it should come as no surprise that in order to prove nonexistence of positive solutions of \eqref{neq} we always have to make a hypothesis on the behavior of $f(t)$ at $t=0$. It turns out that close to zero $f(t)$ should be no worse than a power $t^\sigma$, where $\sigma=\sigma(\lambda, N)$ tends to $2_*=N/(N-2)$ when $\lambda\to1$, and $\sigma$ tends to $1$ when $\lambda\to0$. In addition, we discover that when the ellipticity is too bad (that is, $\lambda$ is too close to zero depending on~$N$), we need to impose a condition on $f(t)$ at $t=\infty$ as well.

\begin{thm}[\cite{AS}]\label{liopuc} Set $\displaystyle\lambda^*=\frac{N-1+\lambda}{N-1-\lambda}$ and suppose that
$$\liminf_{t \to 0}\, t^{-\lambda^*}f(t) > 0.$$
In addition, assume that
\begin{enumerate}
\item[(i)]  if $\lambda=\frac{1}{N-1}$, then $\liminf_{t \to \infty}\, e^{a t} f(t) > 0$ for each $a> 0$; and
\item[(ii)] if $\lambda<\frac{1}{N-1}$, then $\liminf_{t \to \infty} \, t^{|\lambda_*|}f(t) > 0$, where $\displaystyle\lambda_*=\frac{N-1+1/\lambda}{N-1-1/\lambda}$.
\end{enumerate}
Then the inequality \eqref{neq} has no positive weak solution.
\end{thm}

The second hypothesis in Theorem \ref{liopuc} is not very strong. It is needed only when $\lambda<1/(N-1)$ and allows $f(t)$ to decay to zero when $t$ goes to infinity, but no worse than $t^{-|\lambda_*|}$.

Theorem \ref{liopuc} is again optimal, in the sense that we can construct a solution of $\eqref{neq}$, provided we take $f$ to be a model nonlinearity which does not satisfy one of the hypotheses of the theorem (see for instance  \cite{AS0}).

All previous papers on nonexistence for inequalities in non-divergence form concerned the nonlinearity $f(t)=t^q$ (with the exception of \cite{AS0} where we  imposed a more general but still global hypothesis on $f$). A list of references is given in \cite{AS}; we only mention here that it was proved by Cutri and Leoni \cite{CL} that the inequality $-\Puccisub (D^2 u) \geq u^q$ has no positive solutions in the whole space $\rn$ provided $q\in (0, \lambda^*]$. It follows in particular from Theorem \ref{liopuc} that this inequality has no solutions even in any exterior domain of $\rn$, for the larger range $q\in (-\infty,\lambda^*]$ if $\lambda(N-1)\ge 1$, and $q\in [\lambda_*,\lambda^*]$ if $\lambda(N-1)<1$. Of course, Theorem \ref{liopuc} goes much further, by showing that only the behavior of $f(t)$ close to $t=0$ and $t=\infty$ matters, and by describing with precision the behavior which may be allowed.

\smallskip

Theorems \ref{lio} and \ref{liopuc} are very particular cases of Corollary 4.2 in \cite{AS}. The proof of this result is based on a new argument which, in addition to yielding new and optimal results on nonexistence, has several advantages for proving these kinds of Liouville theorems. Above all, it is based entirely on very general maximum principle ideas, which renders it applicable to a great variety of elliptic equations and systems, set in various unbounded domains. We have shown in \cite{AS} how our method trivially extends to systems of elliptic inequalities in exterior domains, still giving optimal results for such systems. We also show in \cite{AS} and in our forthcoming work \cite{ASS} that it yields new nonexistence results in conical domains, and explains the  somewhat different phenomena which occur in such domains.

Next,  besides its obvious simplicity,  the argument makes very apparent the interplay between the scaling of the differential inequality and the scaling of any given subsolution of the differential operator. Optimal results are obtained when this subsolution is taken to be the {\it fundamental solution} of the operator. Finally, our method is independent of the nature of the equations considered, in divergence or non-divergence form, or of the nature of their weak solutions, as long as they satisfy a weak comparison principle. It is actually possible to axiomatize the properties of the elliptic operators involved, under which the method can be applied. We refer to \cite{AS} for a discussion; we expect variations of our method to apply to even larger classes of inequalities.
\medskip

In the next section we describe the proof of Theorem \ref{lio}, dividing it into three parts. We start by giving a list of its main ingredients, then prove some simple particular cases of the theorem which require only subsets of these ingredients, and finally we expose the full proof.

\section{Proof of Theorem~\ref{lio}}

In this section we give the proof of Theorem \ref{lio}. The proof of Theorem \ref{liopuc} is practically the same, see the end of this section. To fix ideas, in the sequel we assume $G\subset B_{1}$ (for each $r>0$ we denote with $B_r$ the ball of radius $r$) and set $\Phi_p(x)=|x|^{(p-N)/(p-1)}$ if $p\not=N$, and $\Phi_p(x)= \log(3|x|)$ if $p=N$.
\bigskip

The basic idea of the proof of Theorem \ref{lio} is very simple. The term $f(u)$ on the right side of \eqref{peq} forces a hypothetical supersolution $u$ of \eqref{peq} to be small. This is because, as for example in the case $f$ is superlinear near $t=0$ and $p < N$, if $u$ were not small then the right-side would be too big for the left side of \eqref{peq}. On the other hand, by the comparison principle, the fundamental solution provides a lower bound for $u(x)$ for large $|x|$. This can be seen by ``sliding" the fundamental solution $\Phi_p$ underneath $u$. These two forces are obviously in conflict, and we would like to understand when this conflict is fatal to the existence of $u$.

\medskip

\subsection{The ingredients of the proof} The key tool we use in estimating $u$  is the following growth lemma, which is a quantitative version of the strong maximum principle. For an easy proof we refer to \cite[Theorem 3.3]{AS}.

\begin{lem}\label{growth}
Assume $h \in L^\infty(B_3\setminus B_{1/2})$ is nonnegative, and $u\geq 0$ satisfies
\begin{equation*}
-\Delta_p u \geq h(x)  \quad \mbox{in} \ B_3\setminus B_{1/2}.
\end{equation*}
(a)\ For each $A\subset B_3\setminus B_{1/2}$ there exists a constant $c_0>0$ depending only on $N$ and $|A|$, such that
$$
\inf_{B_2\setminus B_1}u \geq c_0 \left( \inf_A h\right)^{1/(p-1)}.
$$
(b) Suppose in addition that $u\ge k\Phi_p$ in $B_3\setminus B_{1/2}$ for some $k>0$. Then we have the estimate
\begin{equation*}
 \inf_{B_2\setminus B_1}(u- k\Phi_p)\geq c_0 \left( \inf_A h\right)^{1/(p-1)}.
\end{equation*}

\end{lem}

The fundamental solution $\Phi_p$ and its opposite $-\Phi_p$ give  bounds on the decay (or growth) of any positive $p$-superharmonic function in an exterior domain. This is summarized in the next two lemmas which cause the hypothesis of Theorem~\ref{lio} to break into the different cases it does.
These lemmas are known, though not so often used; we refer to \cite[Lemma 3.7]{AS} for simple proofs based on the comparison principle.

For each $r\ge 1$, we have $\Phi_p>0$  in $\R^N\setminus B_r$,  and so  we may define
\begin{equation} \label{defrrho}
m(r) : = \inf_{B_{2r} \setminus B_r} u, \qquad \rho(r) := \inf_{B_{2r} \setminus B_r} \frac{u}{\Phi_p}, 
\end{equation}

 \begin{lem} \label{mbounds1} Assume that  $u\geq 0$ is $p$-superharmonic in $\rn\setminus B_{1}$. Then
\begin{equation*}
\begin{cases}
r\mapsto \rho(r) \ \mbox{is nondecreasing on} \ [1,\infty), & \mbox{if } \; p < N,\\
r\mapsto \rho(r) \ \mbox{is bounded on} \ [2,\infty), & \mbox{if } \; p \geq N.
\end{cases}
\end{equation*}
In particular $m(r)\ge c r^{(p-N)/(p-1)}$ if $p<N$, $m(r)\le C \log r$ if $p=N$, and $m(r)\le C r^{(p-N)/(p-1)}$ if $p> N$, where $c,C>0$ are constants independent of $r$.
\end{lem}

 Lemma \ref{mbounds1} is proved by sliding underneath $u$  functions of the type $A\Phi_p + B$, for suitably chosen constants $A>0,B\in \mathbb{R}$.

\begin{lem} \label{mbounds2} Assume that  $u\geq 0$ is $p$-superharmonic in $\rn\setminus B_{1}$. Then
\begin{equation*}
\begin{cases}
r\mapsto m(r) \ \mbox{is bounded on} \ [2,\infty), & \mbox{if } \; p < N,\\
r\mapsto m(r) \ \mbox{is nondecreasing on} \ (1,\infty), & \mbox{if } \; p \geq N.
\end{cases}
\end{equation*}
\end{lem}

To prove Lemma \ref{mbounds2} we compare $u$ from below with functions of the type $A(-\Phi_p) + B$, for suitable $A>0,B\in \mathbb{R}$.
\medskip

Even if we will not use it here, for completeness we recall the following simple property, which may be combined with the above lemma.
\begin{prop}\label{propo} Assume that  $u\geq 0$ is $p$-superharmonic in $\rn$ or that $p\leq N$ and $u\geq 0$ is $p$-superharmonic in $\rn\setminus \{0\}$. Then $m(r)$ is nonincreasing in $r>0$.
\end{prop}

\medskip

In order to get the optimal local behavior of the nonlinearity $f$, we will also need the following measure theoretic estimate.

\begin{lem}\label{vwhlap}
For every $0 < \nu < 1$, there exists a constant $\bar C = \bar C(N,\nu)>1$ such that for any positive $p$-superharmonic function $u$ in $B_3 \setminus \bar B_{1/2}$ and any $x_0 \in B_2 \setminus B_1$, we have
\begin{equation*}
\left| \left\{ u \leq \bar C u(x_0) \right\} \cap (B_2\setminus B_1) \right| \geq \nu \left|B_2 \setminus B_1\right|.
\end{equation*}
\end{lem}

We remark that Lemma \ref{vwhlap} is easily seen to be weaker than the weak Harnack inequality proved for example in \cite{Se,T}. See also \cite[Remark 3.6]{AS}.

\medskip

We will next show how these ingredients can be combined to yield nonexistence results. Before giving the full proof of Theorem \ref{lio} which requires all lemmas above, we are going to prove several particular cases which are particularly simple (but still more general than what is usually encountered in the literature) and which need only a subset of these lemmas. We do this in order to, on one hand, better highlight the main points in the proofs, and on the other hand, facilitate eventual extensions of our method to situations in which not all of the above lemmas are available.

\subsection{Some particular cases of Theorem \ref{lio}} Let us first assume $p<N$. In order to simplify the following proofs we will strengthen (i) and assume that
\begin{equation}\label{lim1}
\lim_{t \to 0}\, t^{-p_*}f(t) =\infty.
\end{equation}

\noindent \textbf{1.}
First we are going to show that only Lemmas \ref{growth} and \ref{mbounds1} are sufficient to prove that \eqref{peq} has no positive solutions provided \begin{equation}\label{cond1}
\liminf_{t \to \infty}\, \frac{f(t)}{t^{p-1}} >0.
\end{equation}

\begin{proof} For $r\ge  2$ we define the rescaled function $u_r (x) := u(rx)$ and observe that \eqref{peq} may be written in terms of $u_r$ as
\begin{equation*}
-\Delta_p u_r \geq r^p f(u_r) \quad \mbox{in} \ \R^N \setminus B_{1/r}\supset B_3\setminus B_{1/2}.
\end{equation*}
Since $m(r)=\inf_{B_2\setminus B_1} u_r$ we immediately obtain
\begin{equation*}
-\Delta_p u_r \geq r^p \left(\min_{t\in [m(r),\infty)}f(t)\right)\chi_{B_2\setminus B_1} \quad \mbox{in } B_3\setminus B_{1/2},
\end{equation*}
where $\chi_Z$ denotes the characteristic function of a set $Z\subset\rn$. Applying Lemma \ref{growth} (a) with $A=B_2\setminus B_1$ we  deduce that
\begin{equation}\label{ineq1}
m(r)^{p-1} \geq c_0 r^p \left(\min_{t\in [m(r),\infty)}f(t)\right),
\end{equation}
for some $c_0>0$ independent of $r$. By \eqref{cond1} the  minimum in the right-hand side of \eqref{ineq1} is attained, say at a point $\bar m(r)\in [m(r),\infty)$. So \eqref{ineq1} implies
\begin{equation}\label{ineq2}
\bar m(r)^{1-p} f(\bar m(r)) \leq c_0^{-1} r^{-p} ,
\end{equation}
Sending $r\to \infty$ in \eqref{ineq2} we see that \eqref{cond1}, the continuity of $f$ and $f(t)>0$ for $t>0$ imply $\bar m(r)\to 0$ as $r\to \infty$. Hence \eqref{lim1} implies that for sufficiently large $r$ we have \begin{equation}\label{ineq3}
f(\bar m(r))\ge C_r \bar m(r)^{p_*},
 \end{equation}
 where $C_r\to \infty$ as $r\to \infty$.
Recalling the definition of $p_*$ (which in particular implies that $p-1-p_* = -p(p-1)/(N-p)$) and combining \eqref{ineq3} with \eqref{ineq2} and the lower bound on $m(r)$ in Lemma \ref{mbounds1} yields
\begin{equation}\label{ineq4}
C_r r^{-p}\le  C_r m(r)^{p_*-p+1} \le C_r\bar m(r)^{p_*-p+1}  \leq c_0^{-1} r^{-p} ,
\end{equation}
which is a contradiction, since $C_r\to \infty$ as $r\to \infty$.
\end{proof}
\medskip

\noindent \textbf{2.} We now observe that adding Lemma \ref{mbounds2} to the above argument permits us to relax the hypothesis \eqref{cond1} to the following one
\begin{equation}\label{cond2}
\liminf_{t \to \infty}\, f(t)>0.
\end{equation}

\begin{proof} By Lemma \ref{mbounds2} the left-hand side of \eqref{ineq1} is bounded above, so  the minimum in the right-hand side  tends to zero as $r\to \infty$. Again by \eqref{cond2} and $f(t)>0$ for $t>0$ this minimum is attained and we get $\bar m(r)\to 0$ as $r\to \infty$, so \eqref{ineq3}-\eqref{ineq4} hold, yielding the same contradiction as above. \end{proof}
 \medskip

 Suppose now that $p\ge N$. Note we obtained \eqref{ineq1} by using only Lemma \ref{growth}, so this inequality is independent of the value of $p>1$. In the case $p\ge N$ Lemma \ref{mbounds2} implies that $m(r)$ is bounded below by a positive constant for large $r$. Hence if we assume \eqref{cond2} we see that the minimum in the right-hand side of \eqref{ineq1} is bounded below by a positive constant, by the continuity and positivity of $f$.
 Then \eqref{ineq1} implies that $m(r)\ge cr^{p/(p-1)}$ which is a contradiction with the upper bound in Lemma \ref{mbounds1}, for all sufficiently large $r$. We arrive at the same contradiction if we use only Lemma \ref{growth} and Lemma \ref{mbounds1}, but assume in addition to \eqref{cond2} that $\liminf_{t \to 0}\, \frac{f(t)}{t^{p-1}} >0$.
 \smallskip

We remark that in the  simple proofs above we used only the trivial particular case of Lemma \ref{growth} (a) with $A=B_2\setminus B_1$.
\medskip

Let us now move to the hard part of the proof of Theorem 1, the removal of any condition on $f$ at infinity if $p<N$, resp. allowing $f$ to decay at infinity if $p\ge N$. It is here that we need the very weak Harnack inequality (Lemma \ref{vwhlap}) as well as the growth lemma (Lemma \ref{growth}) in its full strength.

\subsection{Proof of Theorem~\ref{lio}}
Let us suppose that $u> 0$ is a supersolution of \eqref{peq} in $\R^N \setminus B_{1}$. We consider first the case that $p < N$.

\medskip

For $r\ge 2$ we again define the rescaled function $u_r (x) := u(rx)$ and recall  that \eqref{peq} may be written in terms of $u_r$ as
\begin{equation*}
-\Delta_p u_r \geq r^p f(u_r) \quad \mbox{in} \ \R^N \setminus B_{1/r}\supset B_3\setminus B_{1/2}.
\end{equation*}
Denote
\begin{equation*}
A_r:=\{x\in B_2\setminus B_1: m(r)\le u_r(x)\le \bar C m(r)\},
\end{equation*}
 where $\bar C = \bar C(N,\frac{1}{2})> 1$ is as in Lemma \ref{vwhlap}, which yields that
\begin{equation*}
|A_r|\ge (1/2)|B_2\setminus B_1|=c(N)>0.
\end{equation*}
Lemma \ref{growth} (a) with $A=A_r$ and $h(x):= r^pf(u_r(x))\chi_{A_r}(x)$ yields the estimate
\begin{align}\label{limmins1}
m(r)^{p-1} \geq  c_0 r^p \min \left\{ f(t) : m(r) \leq t \leq \bar C m(r)\right\} \quad \mbox{for each} \ r\geq 2,
\end{align}
and some $c_0>0$ independent of $r$.
According to Lemma \ref{mbounds2}, the quantity $m(r)$ is bounded above for $r>2$, and so we have
\begin{equation} \label{limmins}
\min_{\left[ m(r),\bar C m(r)\right]} f\le C r^{-p}m(r)^{p-1}\le Cr^{-p} \to 0\quad\mbox{ as }\; r\to \infty .
\end{equation}
 We deduce from $m(r)\le C$ that the interval over which this minimum is taken is bounded above. Since $f$ is continuous and positive on $(0,\infty)$, $m(r_n)\ge c>0$ for some sequence $r_n\to\infty$ is clearly in contradiction with \eqref{limmins}. Hence
\begin{equation}\label{tooo}
 m(r) \to 0 \quad \mbox{as} \ r \to \infty.
 \end{equation}
 We remark that if instead of (i) we assumed the stronger hypothesis \eqref{lim1}, at this stage  we can deduce the inequalities \eqref{ineq2}-\eqref{ineq4} (with $\bar m(r)\in [m(r),\bar C m(r)]$) and hence a contradiction.

Let us continue with the proof of the Theorem assuming only (i). This assumption and \eqref{limmins1} imply that for $r>2$ sufficiently large,
\begin{equation}\label{limmins2}
(m(r))^{p-1} \geq c r^{p} m(r)^{p_*}.
\end{equation}
Since $p-1-p_* = -p(p-1)/(N-p)$, we can rearrange this inequality as
\begin{equation} \label{mupbnd}
m(r) \leq C r^{-\alpha^*} \quad \mbox{for every sufficiently large} \ r\geq 2.
\end{equation}
where, as before, $\alpha^*=(N-p)/(p-1)$. Note \eqref{mupbnd} and the fact that $\Phi_p$ is $-\alpha^*$-homogeneous imply that the quantity $\rho(r)$ is bounded above for $r>2$.

We recall that, according to Lemma \ref{mbounds1}, for some $c> 0$ we also have
\begin{equation} \label{mlwbnd}
m(r) \geq c r^{-\alpha^*} \quad \mbox{for every} \ r > 2.
\end{equation}

Observe that by Lemma~\ref{mbounds1} $\rho(r)$ is nondecreasing in $r$, that is $\rho(r) = \displaystyle \inf_{\rn \setminus B_r} \frac{u}{\Phi_p}$.
We are going to apply the second part of Lemma~\ref{growth} to the function
\begin{equation*}
v_r(x) :=  r^{\alpha^*}u(rx),
\end{equation*}
which satisfies $v_r\ge \rho(r/2)\, r^{\alpha^*}\Phi_p(rx)=\rho(r/2)\, \Phi_p(x)$ on $\rn\setminus B_{1/2}$.

By \eqref{mlwbnd} we have $v_r\ge c>0$ on $B_2\setminus B_1$, where $c>0$ does not depend on $r$. It also follows from the above argument (inequalities \eqref{limmins1}-\eqref{mlwbnd}) that
$$
-\Delta_p u \ge \frac{f(u)}{u^{p_*}} u^{p_*}\ge c u^{p_*}\quad \mbox{on }\ \frac{1}{r}A_r.
$$
By the scaling invariance and the choice of $p_*$ and $\alpha^*$ (recall the discussion following \eqref{mod}) we see that $v_r$ then  satisfies
$$
-\Delta_p v_r \ge  c v_r^{p_*}\quad \mbox{on }\ A_r,
$$
hence
$$
-\Delta_p v_r \ge c\chi_{A_r} \quad \mbox{on }\ B_3\setminus B_{1/2}.
$$
By Lemma~\ref{growth} (b) this implies
$$
v_r\ge \rho(r/2)\Phi_p + c\ge (\rho(r/2) + c)\Phi_p\quad \mbox{on }\ B_2\setminus B_{1},
$$
where $c>0$ does not depend on $r$. The definition of $\rho(r)$ and the last inequality yield
$$\rho(r) \geq \rho(r/2) + c_1$$ for each $r>2$, and therefore $\lim_{r\to \infty}\rho(r) = \infty$, which contradicts our inequality \eqref{mupbnd}.

\medskip

We next consider the case $p \geq N$. As before, we arrive at the inequality \eqref{limmins1}. Now applying Lemma~\ref{mbounds2}, we see that $m(r)$ is bounded away from $0$, in contrast to the previous case $p < N$ above. Since $f(t)$ is continuous and positive for $t>0$ it follows from the first inequality in \eqref{limmins} that 
\begin{equation}\label{toinf}
 m(r) \to \infty \quad \mbox{as} \ r \to \infty.
 \end{equation}  We split the remainder of the argument into the cases $p=N$ and $p> N$.

In the case $p=N$, we deduce from \eqref{limmins1} and (ii) that for every $a> 0$ there exist $\ep=\ep(a)>0$ small enough and $r_0> 1$ large enough that
\begin{equation*}
 m(r)^{p-1}e^{a \bar C m(r)} \geq c \ep r^p.
\end{equation*}
for $r>r_0$. This inequality trivially implies that for each $A>0$ there exists $r_1>r_0$ large enough that
\begin{equation*}
m(r) \geq A \log r
\end{equation*}
for $r>r_1$, which contradicts the upper bound in Lemma~\ref{mbounds1}.

Finally, in the case that $p > N$, from \eqref{limmins1} and (iii) we obtain
\begin{equation*}
m(r) \geq c r^{(p-N)/(p-1)}
\end{equation*}
for some constant $c> 0$ and large enough $r$. Lemma~\ref{mbounds1}  gives the reverse inequality, so that for large $r$ we have the two-sided bound
\begin{equation} \label{fgh}
cr^{-\alpha^*} \leq m(r) \leq C r^{-\alpha^*},
\end{equation}
that is, $0<c\le \rho(r)\le C$ for all large $r$. We set again  $v_r(x) :=  r^{\alpha^*}u(rx)$, note that $0<c\le v_r\le C$ on the set $A_r$, and argue exactly like in the case $p<N$ to deduce $\rho(r) \to +\infty$ as $r\to \infty$, in contradiction to \eqref{fgh}. The proof is complete.\hfill $\Box$
\bigskip

\begin{remark}
In the case when $p\not=N$ the function $(\Phi_p)^\tau$ is a solution of $-\Delta_p u = u^q$ in $\rn\setminus\{0\}$, for a suitable chosen power $\tau$, provided (i) or (iii) is not satisfied. By using a cut-off like argument it is possible to modify these functions to obtain solutions of $-\Delta_p u\ge u^q$ in $\rn$ (see for instance the argument on page 11 in \cite{AS0}).
If $p=N$ it is also easy to obtain solutions of $-\Delta_p u= e^{au}$, $a>0$, in exterior domains. For instance
\begin{equation*}
u(x) := \textstyle\frac{2}{a}\left( \log|x| + \log\left( \log|x| \right) \right)
\end{equation*}
is a positive solution of the equation
$
-\Delta u = \textstyle\frac{2}{a} e^{-a u}$ in $\mathbb{R}^2\setminus B_{3}$.
Note that for $p=N$ every positive solution of $-\Delta_p u\ge 0$ in $\rn\setminus\{0\}$ is constant (by Lemma \ref{mbounds2} and Proposition \ref{propo}).

These remarks attest to the sharpness of Theorem 1.
\end{remark}\medskip

Finally, we observe that the same argument leads to the proof of Theorem \ref{liopuc}. We only need to replace the fundamental solution $|x|^{(N-p)/(p-1)}$ by $|x|^{\lambda^{-1}(N-1)-1}$, and its opposite $-|x|^{(N-p)/(p-1)}$ by $-|x|^{\lambda(N-1)-1}$ (of course some care is needed about  the different homogeneities of the functions $|x|^{\lambda^{-1}(N-1)-1}$, $-|x|^{\lambda(N-1)-1}$, which both solve $\Puccisub(u)= 0$ in the punctured space $\rn\setminus\{0\}$). Theorem \ref{liopuc} is also sharp, as adequate powers of these functions show.
\bigskip

\noindent {\bf Acknowledgement}. The first
author was supported in part by NSF Grant DMS-1004645.

\bibliographystyle{plain}
\bibliography{tinyliouv}

\begin{thebibliography}{10}

\bibitem{AS0}
Scott~N. Armstrong and Boyan Sirakov.
\newblock Liouville results for fully nonlinear elliptic equations with power
  growth nonlinearities.
\newblock {\em Ann. Sc. Norm. Super. Pisa Cl. Sci. (5)}, to appear.

\bibitem{AS}
Scott~N. Armstrong and Boyan Sirakov.
\newblock Nonexistence of positive supersolutions of elliptic equations via the
  maximum principle.
\newblock {\em Comm. Partial Differential Equations}, to appear.

\bibitem{ASS}
Scott~N. Armstrong, Boyan Sirakov, and Charles Smart.
\newblock Singular solutions of fully nonlinear elliptic equations in cones.
\newblock In preparation.

\bibitem{BP}
M.-F. Bidaut-V{\'e}ron and S.~Pohozaev.
\newblock Nonexistence results and estimates for some nonlinear elliptic
  problems.
\newblock {\em J. Anal. Math.}, 84:1--49, 2001.

\bibitem{B}
Marie-Fran{\c{c}}oise Bidaut-V{\'e}ron.
\newblock Local and global behavior of solutions of quasilinear equations of
  {E}mden-{F}owler type.
\newblock {\em Arch. Rational Mech. Anal.}, 107(4):293--324, 1989.

\bibitem{CL}
Alessandra Cutr{\`{\i}} and Fabiana Leoni.
\newblock On the {L}iouville property for fully nonlinear equations.
\newblock {\em Ann. Inst. H. Poincar\'e Anal. Non Lin\'eaire}, 17(2):219--245,
  2000.

\bibitem{AM}
Lorenzo D'Ambrosio and Enzo Mitidieri.
\newblock A priori estimates, positivity results, and nonexistence theorems for
  quasilinear degenerate elliptic inequalities.
\newblock {\em Advances in Math.}, 224(3):967--1020, 2010.

\bibitem{G}
Basilis Gidas.
\newblock Symmetry properties and isolated singularities of positive solutions
  of nonlinear elliptic equations.
\newblock In {\em Nonlinear partial differential equations in engineering and
  applied science}, volume~54 of {\em Lecture Notes in Pure and Appl. Math.},
  pages 255--273. Dekker, New York, 1980.

\bibitem{KLS}
V.~Kondratiev, V.~Liskevich, and Z.~Sobol.
\newblock Positive solutions to semi-linear and quasi-linear elliptic equations
  on unbounded domains.
\newblock In {\em Handbook of differential equations: stationary partial
  differential equations}, volume~6, pages 255--273. Elsevier, 2008.

\bibitem{MP}
{E}. Mitidieri and S.~I. Pokhozhaev.
\newblock A priori estimates and the absence of solutions of nonlinear partial
  differential equations and inequalities.
\newblock {\em Tr. Mat. Inst. Steklova}, 234:1--384, 2001.

\bibitem{NS}
Wei-Ming Ni and James Serrin.
\newblock Nonexistence theorems for quasilinear partial differential equations.
\newblock In {\em Proceedings of the conference commemorating the 1st
  centennial of the {C}ircolo {M}atematico di {P}alermo ({P}alermo, 1984)},
  number~8, pages 171--185, 1985.

\bibitem{Se}
James Serrin.
\newblock Isolated singularities of solutions of quasi-linear equations.
\newblock {\em Acta Math.}, 113:219--240, 1965.

\bibitem{SZ}
James Serrin and Henghui Zou.
\newblock Cauchy-{L}iouville and universal boundedness theorems for quasilinear
  elliptic equations and inequalities.
\newblock {\em Acta Math.}, 189(1):79--142, 2002.

\bibitem{T}
N.~S. Trudinger.
\newblock On {H}arnack type inequalities and their application to quasilinear
  elliptic equations.
\newblock {\em Comm. Pure Appl. Math.}, 20:721--747, 1967.

\bibitem{V}
Laurent V{\'e}ron.
\newblock {\em Singularities of solutions of second order quasilinear
  equations}, volume 353 of {\em Pitman Research Notes in Mathematics Series}.
\newblock Longman, Harlow, 1996.

\end{thebibliography}

\end{document}